\documentclass[11pt,twoside]{amsart}

\usepackage{color}
\usepackage{graphicx}

\usepackage{amssymb}
\usepackage{latexsym}

\makeatletter

\@addtoreset{equation}{section} \makeatother
\newtheorem{theorem}{Theorem}[section]
\newtheorem{remark}[theorem]{Remark}
\newtheorem{lemma}[theorem]{Lemma}
\newtheorem{proposition}[theorem]{Proposition}
\newtheorem{corollary}[theorem]{Corollary}
\newtheorem{definition}[theorem]{Definition}
\newtheorem{example}[theorem]{Example}

\def\1{\mathbf{1}}
\def\:{\lrcorner}
\def\#{\sharp}

\def\d{\delta}
\def\e{\epsilon}

\def\o{\circ}
\def\s{\sigma}

\def\qed{\ensuremath{\quad\Box\quad}}

\def\inv#1{\raise.1em\hbox to 0pt{$^{-1}$\hss}_{#1}\;}

\def\grad{{\rm grad}}
\def\v{\noindent}

\newcommand{\bean}{\begin{eqnarray*}}
\newcommand{\eean}{\end{eqnarray*}}
\newcommand{\benu}{\begin{enumerate}}
\newcommand{\eenu}{\end{enumerate}}
\newcommand{\eea}{\end{eqnarray}}
\newcommand{\bea}{\begin{eqnarray}}

\newtheorem{Theorem}{Theorem}

\newtheorem{Definition}{Definition}

\newcommand{\be}{\begin{equation}}
\newcommand{\ee}{\end{equation}}

\newcommand{\N}{{\mathbb N}}
\newcommand{\Z}{{\mathbb Z}}
\newcommand{\R}{{\mathbb R}}

\newcommand{\gr}{{\rm grad}}

\newcommand{\ben}{\begin{enumerate}}
\newcommand{\een}{\end{enumerate}}
\newcommand{\bit}{\begin{itemize}}
\newcommand{\eit}{\end{itemize}}
\newcommand{\edoc}{\end{document}}

\newcommand{\bdefi}{\begin{definition}}
\newcommand{\btheo}{\begin{theorem}}
\newcommand{\bprop}{\begin{proposition}}
\newcommand{\brema}{\begin{remark}}
\newcommand{\bcoro}{\begin{corollary}}
\newcommand{\blemm}{\begin{lemma}}
\newcommand{\bexam}{\begin{example}}

\newcommand{\edefi}{\end{definition}}
\newcommand{\etheo}{\end{theorem}}
\newcommand{\eprop}{\end{proposition}}
\newcommand{\erema}{\end{remark}}
\newcommand{\ecoro}{\end{corollary}}
\newcommand{\elemm}{\end{lemma}}
\newcommand{\eexam}{\end{example}}

\begin{document}

\setcounter{page}{1}

\title{Special temporal functions on globally hyperbolic manifolds}

\author[Olaf M\"uller]{Olaf M\"uller}
\address{Fakult\"at f\"ur Mathematik, Universit\"at Regensburg, Universit\"atsstra\ss e 31, 93053 Regensburg, Germany}
\email{olaf.mueller@mathematik.uni-regensburg.de}

\date{\today}

\parindent=5mm

\maketitle

\begin{abstract}
\v In this article, existence results concerning temporal functions with additional properties on a globally hyperbolic manifold are obtained. These properties are certain bounds on geometric quantities as lapse and shift. The results are linked to completeness properties and the existence of closed isometric embeddings in Minkowski spaces.

\end{abstract}

\section{Introduction}

\v In the research on classical field theory on Lorentzian manifolds, the most appropriate geometric category for classical field theory turned out to be the one of oriented globally hyperbolic manifolds and their causal isometric oriented and time-oriented embeddings. On one hand, this is due to its relatively easy and invariant definition as the category of manifolds with compact causal diamonds and without closed causal curves (for this definition weaker than the usual one see \cite{BeSa}), on the other hand to the strong statements about the well-posedness of initial-value problems of Laplace-type or normally-hyperbolic operators on them, i.e., of operators whose symbols coincide with the Lorentzian metric tensored with the identity in the configuration bundle of the respective field theory (\cite{BGP}, for extensions of results of this kind to larger classes of field theories cf. \cite{vW}, e.g.). An important property of a globally hyperbolic manifold is that it admits an orthogonal foliation by Cauchy surfaces (\cite{BeSa1}, \cite{BeSa2}), that means, it is isometric to $(\R \times S, G = - f^2 dt^2 + g_t)$, where $t : \R \times S \rightarrow \R$ is the projection on the first factor, $f : \R \times S \rightarrow (0, \infty)$ is a smooth function, and $g_t $ is a smooth family of Riemannian metrics on $S$. Moreover, all level sets of the function $t$ are Cauchy surfaces. Now, in the initial value formulation of general relativity, the initial data are two 2-tensors $g_0, W$ on a 3-manifold representing the induced Riemannian metric on the Cauchy surface $\Sigma$ with normal vector field $n$ and the Weingarten tensor of the Cauchy surface, respectively, subject to the constraint equations (where $T$ is the stress-energy tensor)

$$  R^{\Sigma} + (tr_{g_0} W)^2 - \vert \vert W \vert \vert_{g_0}^2 = 16 \pi T(n,n), \qquad  tr_{1,3} (\nabla W) - tr_{1,2} (\nabla W) = 8 \pi T(n, \cdot)   $$

\v Thus bounds on these quantities on Cauchy surfaces are of natural interest.

\section{Definitions and relation to embeddings}

\v A function $t: M \rightarrow \R$ on $M$ is called {\bf time function} iff it is non-decreasing along every future timelike curve (if $t$ is $C^1$ this is equivalent to the requirement that its gradient ${\rm grad} t $ be past-directed causal or zero at every point). It is called {\bf temporal function} iff it is $C^1$ and its gradient ${\rm grad} t $ is past-directed timelike at every point. Let us introduce the notions $ P_a (t) := t^{-1} ((- \infty, a)) $ and $ F_a (t) := t^{-1}( (a, \infty)$ as well as $S_a (t) := t^{-1} (\{ a \})$. If, for a time function $t$, the function $ t \o c$ is surjective onto $\R $ for any inextendible causal curve, we call $t$ {\bf Cauchy} (as this is equivalent to require that the preimage of every real number is a Cauchy surface). Now a smooth Cauchy temporal function $t$ gives rise to an isometry between $M$ and $\R \times N$ as above, and vice versa. For $g = - f^2 dt^2 + g_t$ (here we have $f^{-2} = g ({\rm grad}^g t, {\rm grad}^g t)$), the normal vector field at a level set of $t$ is $n= f^{-1} \partial_t$, and the Weingarten tensor is easily computed as $W(X,Y) = g(\nabla_X n, Y) = - f {\rm Hess}(t) (X,Y)$. Therefore one should ask for special foliations in which these quantities are bounded. The following definition systematizes foliations with different bounds. As a preparation, for some distribution $E$ on $M$, we define the $E$-flip metric $g^E$ by $g^E \vert_{E^{\perp}} := g  \vert_{E^{\perp}}$ and  $g^E \vert_{E} := - g  \vert_{E}$ (this is well-defined as $E \cap E^{\perp}$ consists of lightlike vectors). For a vector field $V$, we write $g^V:=  g^{\R \cdot V }$. Of course, for $D$ timelike, we have $g^D > \vert g \vert $. For a temporal function $t$, the $t$-flip metric $g_+^t$ is defined as $g_+^t := g^{{\rm grad} t}$.

\begin{Definition}
\begin{enumerate}
\item{A smooth function $s$ on $M$ is called {\bf steep in $A$}, for some $A \subset M$, if there is a negative number $h$ such that $g({\rm grad} s, {\rm grad} s) < h $ in $A$. If $s$ is steep in every level surface and if the bound of the $C^k$ norm depends continuously on the level surface then $s$ is called {\bf semi-steep} or an {\bf s function}. It is called {\bf S} iff it is steep in $M$.} 
\item{For any $n \in \N $ and any Riemannian metric $G$ on $M$, a smooth function $s$ on $M$ is called {\bf $(n,G,A)$-mild} iff $\vert \vert {\rm grad} s \vert \vert_{C^n (A,G)} < \infty $ \footnote{The seminorm $\vert \vert  \cdot \vert \vert_{C^k (A,G)}$ applied to a vector field is defined in the usual way by $\vert \vert u \vert \vert_{C^k (A,G)}  := \sum_{i=1}^k {\rm sup} \{ \vert \vert \nabla^{(k)} (V_1, ... , V_k) u (x) \vert \vert_G : x \in A, V_i \in T_x M, G(V_i, V_i) = 1  \} $.}. Omission of entries is defined as insertion of the default arguments $n=1 $, $G= g_+^s$, $A=M$. If a function is mild in every level surface and if the bound of the $C^k$ norm depends continuously on the level surface then $s$ is called {\bf semimild} or an {\bf m function}. A function $s$ which is mild in $M$ is called {\bf mild} or {\bf M function}.}   
\item{A temporal function is called {\bf C} iff its level sets are complete Riemannian manifolds if equipped with the induced Riemannian metric.}
\end{enumerate}
\end{Definition}

\v Obviously, a temporal function is steep if and only if the function $f$ in the corresponding decomposition as above is bounded. As an example for mildness, a function $s$ is $1$-mild if and only if the $g$-length of its gradient is bounded in $ M  $ and if for its Hessian $H$ we have $\vert H(v,v) \vert < C \cdot g_+^s(v,v)$ on $M$ for some $C >0$. Of course, every Cauchy temporal function on a globally hyperbolic manifold whose Cauchy surfaces are compact is trivially $(n,G)$-semimild for each $n$ and $G$, in particular an $s$ function. By convention, accumulation of letters to combinations like SMC attributed to a temporal function means that it has all the respective properties, and attributed to a manifold $M$ it indicates that $M$ has {\em one Cauchy temporal function with all the required properties}, in this case such that the level surfaces are Cauchy, the Hessian is controlled as above, and the shift is bounded from infinity. This is a subtle distinction, as there are examples of globally hyperbolic manifolds with two Cauchy temporal functions one of which has complete levelsets while the other hasn't: Consider Geroch's example in \cite{BEE}, p. 204, and invert the sign of the two-dimensional metric, in other words, consider a Lorentzian metric $g$ on $\R^2$ conformally equivalent to the standard Minkowski metric with conformal factor $ \phi$ with the properties

\begin{enumerate}
\item{$\phi (x_0, x_1) = 1$ for $ \vert x_0 \vert \geq 1 $,}
\item{$\phi (x_0, x_1) = \phi (-x_0, x_1) $ for all $(x_0, x_1 ) \in \R^2$,}
\item{There is a real number $C$ with $ \phi(0, x_1) \leq C \cdot x_1^{-4}$ for all $x_1 \in \R$.}
\end{enumerate}

\v Then, as conformally equivalent to standard Minkowski space, $M:= (\R^2, g)$ is globally hyperbolic, and as every future-pointing geodesic $c$ hitting the strip $ (-1,1 ) \times \R $ in a point $ c(0)$ leaves it forever after crossing the compact subset $J^+ (c(0)) \cap J^- (\{ 1\} \times \R)$. Therefore $M$ is causally geodesically complete. The Cauchy temporal function $x_0$ has the incomplete levelset $ x_0^{-1} (0) $, but as soon as we consider a Cauchy temporal function $T$ given as the Lorentzian scalar product with any future timelike vector not collinear to $\frac{\partial}{\partial x_0}  $, then its level surfaces are all complete. As a second example, an article of Candela, Flores and S\'anchez (\cite{mC}, Sect. 6.2) shows that even usual Minkowski space has non-C Cauchy temporal functions.

\bigskip

\v By rescaling from a fixed level surface we can construct an $M$ Cauchy temporal function from an $m$ Cauchy temporal function in an obvious way, and the property s will be preserved by this rescaling. The same is obviously not true for the property S, and from an s Cauchy temporal function we do not get an S Cauchy temporal function by rescaling as the range may change from all of $\R$ to a proper subset. 

\bigskip

\v In addition to the interest in SMC type decompositions coming from the initial value problem in General relativity, they have also some analytic advantages as long-time existence of minimal surfaces (cf. \cite{oM1}). Note that the definition given in that paper includes that the eigenvalues of $g_t ^{-1} \o \dot g_t$ be bounded on every level set of $t$, but by the bound on $f$ this property follows, as for $g= - f^2 dt^2 + g_t$, and for vector fields $V,X,Y$ invariant under the flow of $\partial_t$ in the corresponding decomposition, we get $  W (X,Y) = -  f \cdot Hess (X,Y) $ and 

\bean
\dot{g}_t (V,V) &=& 2 g(\nabla_t V, V)\\ 
&=& 2 g(\nabla_V \partial_t, V) = 2 g(\nabla_V ((g({\rm grad} t, {\rm grad} t)^{-2} {\rm grad} t , V)\\
&=& 2 g(V(g({\rm grad} t, {\rm grad} t)^{-2}) \cdot {\rm grad} t , V) + 2 g({\rm grad} t, {\rm grad} t)^{-2} g(\nabla_V {\rm grad} t, V)\\
&=& 2 g({\rm grad} t, {\rm grad} t)^{-2} \cdot {\rm Hess} (t) (V,V)\\ 
&=&  - 2 f^2 Hess (V,V)
\eean

\v Apart from the reasons explained above, the additional conditions are also important in the theory of isometric embeddings. In a recent article, Miguel S\'anchez and the author showed the existence of a steep Cauchy temporal function in any globally hyperbolic manifold {\cite{MS}} which implied among other results that it has an isometric embedding in a Minkowski space. Now, for many applications in variational problems, one looks for a {\em closed} embedding of the given manifold into a vector space. It turns out that the existence of a closed embedding into a Minkowski space is equivalent to the SC property:

\begin{Theorem}
\label{embed}
If a Lorentzian manifold of dimension $n$ has an SC Cauchy temporal function then it has a closed isometric embedding into the Minkowski space of dimension $N(n) +2$ where $N(n)$ is the Nash dimension of $n$. If, conversely, $(M,g)$ has a closed isometric embedding into a Minkowski space it has a SC Cauchy temporal function.
\end{Theorem}

\v {\bf Proof.} Let $t$ be the SC Cauchy temporal function, then $2t$ is a SC Cauchy temporal function as well, and $G:= g+ d^2(2t) = g+4dt^2 $ is a Riemannian metric with $\vert \vert X \vert \vert_G > dt (X)$. We want to show its completeness. Let $x_n \in M $ be a $G$-Cauchy sequence, then $\vert t(x_n) - t(x_m) \vert < d^G (x_n, x_m)$. Thus the real numbers $t(x_n)$ form a Cauchy sequence in $\R$ converging to some $T$. By a simple $3 \e$-argument, the sequence $p_T (x_n) $ (where $p_T : M \rightarrow t^{-1} (T)=: S_T$ is the projection onto the level set to $T$ corresponding to the global decomposition) forms a Cauchy sequence in $S_T$ as well: $ d^G (p_T (x_n ), p_T (x_m)) \leq d^G (p_T (x_n), x_n) + d^G(x_n, x_m) + d^G(x_m, p_T(x_m)) $ (and then the first and the third term converge to zero as the $t(x_n)$ converge to $T$ and by the S property). As, by the property C, $S_T$ is complete, this sequence converges to some $p \in S_T$ which is then also a limit of the original sequence $x_n$. This shows completeness of $G$. By the theorem in \cite{o5}, there is a closed isometric embedding $I$ of $(M,G)$ into the Euclidean space of dimension $N(n) +1$. Now add a Lorentzian factor as in \cite{MS}: Define $\tilde{I} (x) := ( 2t (x), I(x))$, which is an isometric embedding for the original, Lorentzian, metric $g$, into $\R^{1, N(n) + 1}$,  then the image of the resulting embedding is a graph over the (closed) image of $I$ and is therefore closed as well. 

\v Conversely, if $(M,g)$ has a closed isometric embedding $I$ into some $\R^{1,n}$, there is a steep temporal function $t$ which is just the restriction of the scalar product with some timelike vector to the submanifold. To show that its level sets are complete, let a Cauchy sequence in a level set $S_a$ be given. Then it is a Cauchy sequence in a Euclidean vector subspace of $\R^{1,n}$ and thus converges in this vector subspace and therefore in $\R^{1,n}$, and moreover, because of its closedness, in the image of the embedding. This limit, finally, is contained in the correct level set because of continuity of $t$. To show that $t $ is Cauchy, let $c$ be an inextendable causal curve in $M$. Then $I \o c$ is a causal curve in $\R^{1,n}$. If $t \o c$ is bounded, then the coordinate function is bounded along $I \o c$, so $I \o c: (-1, 1) \rightarrow \R^{1,n}$ is extendable as a curve in $ \R^{1,n}$ to a larger interval containing $1$. The closedness of $I$ implies that $(I \o c) (1)$ is contained in $i(M)$, therefore $c$ is actually extendable in $M$, contradiction.   \hfill \qed

\section{Results on sm functions}

\v The properties are not independent of each other, as for example the C property for a Cauchy temporal function is a consequence of the properties s and m:

\begin{Theorem}
\label{completelevelsets}
Every level set of an sm Cauchy temporal function $t$ in a time-oriented Lorentzian manifold $(M,g)$ is complete.
\end{Theorem}

\v {\bf Proof.} First, by the existence of a Cauchy temporal function, $(M,g)$ is globally hyperbolic. Recall that the decomposition $M \cong \R \times S$ induced by $t$ carries a metric $-f^2 dt^2 + g_t$ where the function $f$ is bounded below by $2 \d >0$ on every level set by the m property. Because of the S and the m properties there is a constant $D>0$ with $\dot{g}_t (V,V) \leq D g_t(V,V)$ for a vector field $V$ in $S$, cf. the equation before Theorem \ref{embed}. The number $D$ depends continuously on the level set by the continuity requirement in the m property. Consider some piecewise $C^1$ curve $c =: (a, c_2): [0, T) \rightarrow t^{-1}({a}) =:S_a$ parametrized by arc length, then there is a positive $\e$ s.t. for $\theta  \in t^{-1} ((a- 2 \e, a + 2 \e)) $ we have 

$$ \frac{3}{4} g_{\theta} (V,V) \leq g_a (V,V) \leq \frac{5}{4} g_{\theta} (V,V) $$ 

\v for $ t $-invariant $V$ because of the condition on $\dot{g}_t$. This inequality implies that if $c(T - \e \d) = (a, p)$ then for $x:= (a - \e, p)$ we have $c([T - \e \d, T)) \subset I^+ (x)$ as, for every $r < \e \d$, we have $\underline{c}_r: [0, \e \d] \rightarrow M $, $\underline{c}_r (s) := (a - \e + \d^{-1} s, c_2(T - \e \d  + \frac{r}{\e \d} \cdot s) ) $ is a piecewise $C^1$ curve from $x$ to $c(T - \e/2 + r)$ which is future timelike as $k$ is uniformally bounded from $0$. That means that the sequence $ \s : \N \rightarrow M, \s (n) := c(T- 1/n) $, is contained in $J^+ (x) \cap J^- (S_a)$ which is compact and has therefore an accumulation point.  \hfill \qed

\bigskip

\v It would be useful to know how special Sm functions are. The author conjectures that {\em every globally hyperbolic manifold admits a Sm temporal function} (in general, this function cannot be Cauchy, see the section on counterexamples). Attempts to prove this with different techniques encountered surprisingly difficult obstacles. What we, however, {\em can} prove in a comparatively simple manner is the existence of a $(n, H)$-mild function $t$ for each natural $n$, for any Riemannian metric $H$ (whereas the function would be m if it were $(1,g_+^t)$-mild). The precise statement is displayed in the following theorem which could be important for initial value problems as in the introduction:

\begin{Theorem}
\label{haupt}
Let $(M,g)$ be globally hyperbolic, choose a Cauchy surface $S \subset M$. Let $k \in \N$ be given, let $H$ be a Riemannian metric on $M$. Then $(I^+ (S),g \vert_{I^+ (S)})$ has a $(k,H)$-semimild temporal function taking values in $(0,1)$ and approaching $1$ for every $C^0$-inextendible $C^2$ future curve.


\end{Theorem}

\bigskip

\v Before the proof of Theorem \ref{haupt} let us mention a general result closely related to the techniques used in the proof. It concerns the existence of a {\em causal numbering}. Let $(M,g)$ be a spacetime and let $P:= \{ p_i \vert i \in \N \}$ be a discrete sequence of points in $M$. Then an {\bf acausal numbering} of $P$ is a bijection $p$ between $\N$ and $P$ such that $i < j$ implies $p(j) \notin J^+ (p(i))$. In the proof we will need the notion of a set being {\em spatially compact}. A subset $A$ of a globally hyperbolic manifold is called {\bf spatially compact} iff, for any Cauchy surface $S$ of $M$, the sets $A \cap I^+S$, $A \cap I^- ( S)$ are either empty or compact.  

\begin{Theorem}   
Let $(M,g)$ be globally hyperbolic and $S$ be a Cauchy surface of $M$. Then every discrete sequence $P$ in $J^- (S)$ has an acausal numbering.  
\end{Theorem}

\v {\bf Proof.} First we define an exhaustion of $J^- (S) $ by spatially compact future sets $D_i$. We can choose $D_i := D^- (B_i)$ where the $B_i$, e.g. chosen as balls with radius $i$ around a fixed point w.r.t. any complete auxiliary Riemannian metric on $S$, form a compact exhaustion of $S$, and where $D^- $ means the past dependency region (also called 'Cauchy development' in some references), for $A \subset M$ defined by $D^- (A)$ being the subset of $M$ which consists of the points $p$ such that every $C^0$-inextendable timelike future curve starting at $p$ intersects $A$). $D(A)$ is  spatially precompact for compact $A $ but in general not precompact, as it happens in the example of the de Sitter spacetime. We choose any continuous time function $t$ on $J^- (S) $ taking the value $0$ on all of $S$ as described in \cite{BS4}. Now, for finitely many points, obviously there is always an acausal numbering. Using this, we define $P_i : = P \cap A_i$ where $A_i := t^{-1 } ([-i, 0]) \cap D_i$ and $Q_i := P_{i+1} \setminus  P_i$ and paste acausal numberings of $Q_i$ to an acausal numbering of all of $P$. \hfill \qed    

\bigskip








\v {\bf Proof of the Theorem.} Let $ \{ C_n \vert n \in \N  \}$ with $\overline{C_n} \subset {\rm int} C_{n+1}$ for all $n \in \N$ be a compact exhaustion of $S$. Define $A_n := \big( D^- (C_n) \cap t^{-1} ((-n, 0 )) \big) \cup  I^+ (C_n)  $. This is a future subset, in particular a causally convex subset of $M$ and  therefore globally hyperbolic. Now let $\tau_n$ be a smooth Cauchy temporal function of $A_n$ and put $m(n) := {\rm min} \{ \tau_n (x) \vert x \in A_{n-1} \cap I^- (S)  \} $. By adding a constant number to $\tau_n$ if necessary we can satisfy $m (n) >1$. Now define $s_n := \psi \o \tau_n $ with $\psi (r) := e^{-r^{-2}}$ for $r>0$ and $\psi (r) = 0$ for $r<0$. Then $s_n$ (extended by $0$ on $M \setminus A_n$) is a smooth positive function on $M$ with support in $A_n$ which is a temporal function in $A_{n-1} \cap I^- (S) $ while its gradient is always past or zero. As $I^+ (\tau_n^{-1} (\{ m_n\})) \cap I^- (S)$ is precompact in $A_n$, it is $(n,H)$-mild in $I^-(S)$. Then, put

$  s_-:= \sum_{n=1}^{\infty} 2^{-n} \cdot  \min \{ 1, \vert \vert s_n \vert \vert_{C^n (P_n (t))} ^{-1}  \} s_n ,$ 

\v where $P_n(t) = t^{-1} (( - \infty, n))$ as defined in the introduction and the norm on $C^n (P_n(t))$ is the usual $C^n$ norm. The norm is finite as $\s_n$ was supported in $A_n$ and $A_n \cap I^{-1} (S_n)$ is precompact. Finally, a function with the required properties can be constructed by inverting the time direction, constructing the function $s_-$ as above for the inverted time direction and finally defining $ s_+ := 1- \s_i  $  . \hfill \qed 

\bigskip

\section{Relation to b.a.-completeness}

\v The smC properties are also connected with a genuinely Lorentzian notion of completeness called b.a.-completeness. Following the book \cite{BEE}, we define a curve $c : J \rightarrow M$ in a time-oriented Lorentzian manifold $(M,g)$ to be {\bf of bounded acceleration} or {\bf b.a.} iff $ g(\dot{c} (t), \dot{c} (t)) = -1  $ for all $t$ and $ \{  \langle \nabla_t \dot{c} , \nabla_t \dot{c} )  : t \in J \} $ bounded in $\R$, and a time-oriented Lorentzian manifold $(M,g)$ to be {\bf b.a.-complete} if all $C^0$-inextendible $C^2$-curves of bounded acceleration are defined on $\R$, or, equivalently, have infinite length. Of course, b.a. completeness implies timelike geodesic completeness. There are many examples of b.a.-complete spacetimes: 

\begin{itemize} 
\item{Minkowski spaces are b.a.-complete which can be shown by elementary calculation.}
\item{For a b.a.-complete manifold $(M,g)$, any spatially compact perturbation of the metric yields a b.a.-complete manifold again.}
\item{Ehresmann complete fiber bundles with either b.a.-complete fibers and Riemannian base space or with Riemannian fibers and b.a.-complete base space are b.a.-complete (in particular Riemannian products, coverings and warped products).}
\item{compact homogeneous Lorentzian manifolds are b.a.-complete, with very much the same proof as in the book of O'Neill (\cite{O}) where geodesic completeness is shown (recall that $\frac{\nabla}{\partial_t} \dot{c} \perp \dot{c}$!).}
\item{Totally geodesic immersed Lorentzian submanifolds of b.a.-complete manifolds are b.a.-complete.}
\end{itemize}

\v Moreover, the existence of appropriate conformal Killing vector fields also ensures b.a.-completeness as in Prop. 2.1 in \cite{RS} (if restricted to Lorentzian signature) where geodesic completeness was shown. As b.a. completeness implies timelike geodesic completeness but not spacelike geodesic completeness, the result in \cite{RS} and the following one complement each other for Lorentzian manifolds:

\begin{Theorem}
Let $(M,g)$ be a Lorentzian manifold and $X$ a vector field on $M$ with the following properties:

\begin{enumerate}
\item{There is a bounded smooth function $a$ on $M$ such that $L_X g = a \cdot g$ (in particular, $X$ is conformally Killing),}
\item{There is an $\e >0$ such that $g(X,X) < - \e$ (in particular, $X$ is timelike), and}
\item{$g^X$ is complete.}
\end{enumerate}

Then $(M,g) $ is b.a.-complete.
\end{Theorem}

\v {\bf Proof.} Let $c: [0, b) \rightarrow M$, $0 < b < \infty$, be a b.a. curve in $M$. We have to show that $c$ is $C^0$-extendible beyond $b$. Property (3) implies that this is the case if $g^X (\dot{c}, \dot{c}) $ is bounded. Now $g^X (\dot{c}, \dot{c}) = -1 - 2 (g(X, X) )^{-1} g^2 (X, \dot{c})$. As $(g(X, X))^{-1}$ is bounded by Property (2), it remains to show that $g(X, \dot{c}) $ is bounded. We extend $\dot{c} $ to a unit vector field on an open neighborhood of $c([0, b))$ in $M$ (e.g. by using a locally finite covering $C$ of $c([0, b))$ and a partition of unity subordinate to $C$ and the normal exponential map) and calculate

\bean
\frac{d}{dr} (g(X, \dot{c})) &=& g(\nabla_r X, \dot {c}) + g(X, \nabla_r \dot{c})\\
 &=& g(\nabla_X \dot{c}, \dot{c}) + g([\dot{c}, X], \dot{c}) + g(X, \nabla_r \dot{c})\\
&=&  \frac{1}{2} X(g(\dot{c}, \dot{c})) + \frac{1}{2} X(g(\dot{c}, \dot{c})) + (L_X g) (\dot{c}, \dot{c}) +      g(X, \nabla_r \dot{c})\\
&=& ag(\dot{c}, \dot{c}) +    g(X, \nabla_r \dot{c}) = -a +    g(X, \nabla_r \dot{c})\\
\eean

\v We want to estimate the last term by means of $g(X, \dot {c})$. That is indeed possible: we decompose $X$ as $X = X_1 \cdot \dot{c} + X_2 \cdot \nabla_r \dot{c} + W$ with $\nabla_r \dot{c} \perp W \perp \dot{c} $. Thus, with $M:= g(\nabla_r \dot{c}, \nabla_r \dot{c}) < D$ we get

$$ 0 > g(X,X)  = - X_1^2 +  M X_2^2 + g( W, W) \geq M X_2^2 - X_1^2  ,$$

\v thus we have $X_1^2 > M X_2^2 $ which implies $g(X, \dot{c} ) > \sqrt{M } X_2 = M^{-1/2} g(X, \nabla_r \dot{c}) $ (as $g(X, \nabla_r \dot{c}) = g(X_2 \nabla_r \dot{c}, \nabla_r \dot{c}) = X_2 M$) and therefore $g(X, \nabla_r \dot{c}) \leq \sqrt{M} g(X, \dot{c}) $. That means $\vert \frac{d}{dr} (g(X, \dot{c}) ) \vert \leq a + \sqrt{D} g(X, \dot{c})  $ which provides us with the usual exponential estimate for $g(X, \dot{c})$ preventing that it be unbounded within finite time. \hfill \qed

\bigskip

\v Now, if a globally hyperbolic manifold $(M,g)$ is moreover b.a.-complete, then we have the following:

\begin{Theorem}
Let $(M,g)$ be globally hyperbolic and $t$ an SM temporal function on $M$.
\begin{enumerate}
 \item{The integral curves of ${\rm grad t}$, reparametrized by arc length, are of bounded acceleration.} 
\item{If, additionally, $(M,g)$ is b.a.-complete, $t$ is $3$-mild, and the curvature operator of $g$ is $g_+^t$-bounded along every level set of $t$ by a constant continuously depending on the level set, then $t$ is also C and Cauchy (and consequently, $(M,g)$ admits a closed isometric embedding into some Minkowski space).}
\end{enumerate}
\end{Theorem}

\v {\bf Proof.} The mildness property says that for $G:= g_+^t$ we have 

$$  G(\nabla_X^g \grad^g t, \nabla_X ^g \grad^g t  )  \leq D G(X,X) $$

\v for all vectors $X$. Let $ \nabla $ denote the Levi-Cita covariant derivative of $(M,g)$. Then we get, for $\dot{c} (r) = \vert \vert \grad t (c(r)) \vert \vert_g^{-1} \cdot \grad t (c(r))$ and with $A := \grad^g t$,

\bean
g(\nabla_r \dot{c}, \nabla_r \dot{c}) &=& \frac{1}{g(A, A ) ^2} \vert \vert - \frac{1}{g(A, A)} g (\nabla_{A } A , A ) \cdot A + \nabla_A A  \vert \vert_g^2\\
&\leq&  \frac{1}{g(A, A ) ^2} \vert \vert - \frac{1}{g(A, A)} g (\nabla_A A , A ) \cdot A + \nabla_A A  \vert \vert_G^2\\
&\leq&  \frac{1}{g(A, A ) ^2} \cdot \big( g(A , A )^{-2} \cdot \vert \vert  g (\nabla_A A , A ) \cdot A \vert \vert_G^2 + \vert \vert  \nabla_A A   \vert \vert_G^2 \big)\\
\eean

\v where the first inequality stems from the general fact $\vert g \vert < G$ for flip metrics. Now the first term of the sum can be estimated by 

\bean
\vert \vert  g (\nabla_A A , A ) \cdot A  \vert \vert_G^2
\leq G^2 (\nabla_A A, A) \vert \vert A \vert \vert_G^2\\
\leq \vert \vert \nabla_A A \vert \vert_G^2 \cdot \vert \vert A \vert \vert_G^2 \cdot \vert \vert A \vert \vert_G^2
\leq D \vert \vert A \vert \vert_G^6
\eean

\v where the first inequality is due to $g<G$ again and the last one is property M. The second term can be estimated by $D G^2(A, A) = - D g^2(A , A)$ and therefore we have

$$  g(\nabla_r \dot{c}, \nabla_r \dot{c} ) \leq D (\vert \vert \grad^g t \vert \vert_g ^2  + \vert \vert \grad ^g t \vert \vert_g^{-2})  $$

\v which is bounded by the property S and M of $t$. 

\v For the second assertion, assume that there is a non-complete level set of $t$, say, $t^{-1} (\{ a \}) =:S$. Pick a $C^0$-inextendible arclength-parametrized geodesic curve $c: [0,1) \rightarrow S $ (geodesic with respect to the Riemannian metric on $S$). Then we want to construct a b.a.-curve $c_{\e}$ of finite length. To that purpose, we extend the vector field $\dot{c}$ along $c$ to a vector field $X_0$ on a normal neighborhood $U_0$ of the image of $c$ in $S$, and then we extend $X_0$ to a vector field $X$ on an open neighborhood $U$ of the image of $c$ via $\partial_t$-invariance, that is, such that $[\partial_t,X ] = 0  $. Then, for $\e >0$, let a curve $c_{\e}$ be defined by $c_{\e} (0) := c(1 - \e)$, $\dot{c}_{\e} (s) :=  (\partial_t + aX) (c_{\e}) (s)$. Here the function $a$ is determined by the normalization condition $g(\partial_t + aX , \partial_t + a X) = -1$ which is equivalent to 
$  a := \sqrt{\frac{-1-g(\partial_t, \partial_t )}{g(X,X)}} $. Obviously $c_{\e}$ is not $C^0$-extendible beyond $1$ as otherwise $c$ being its projection to $S$ would be. It remains to be shown that $c_{\e}$ is indeed a b.a. curve for $\e$ sufficiently small. Using $\partial_t = g({\rm grad} t, {\rm grad} t)^{-2} \cdot {\rm grad} t $, we compute 

\bean
\nabla_s (\partial_t + a X ) (c_{\e}(s)) &=& \nabla_{\partial_t} \partial_t + \nabla_{\partial_t} aX + \nabla_X \partial_t + \nabla_{aX} aX\\
&=& g({\rm grad} t, {\rm grad} t)^{-4} \nabla_{{\rm grad} t} {\rm grad} t\\
&+& g^{-2} ({\rm grad} t, {\rm grad} t) \cdot   \big( ({\rm grad} t) ( g^{-2}( {\rm grad} t, {\rm grad} t ) )   \big)\\
&+& \partial_t a \cdot \nabla_{\partial_t} X + a \nabla_{\partial_t } X + a \nabla_X \partial_t  + a^2 \nabla_X X + a X(a) \cdot X
\eean

\v First we note that the term $a \nabla_{\partial_t} X  $ can be replaced by $a \nabla_X \partial_t$ by the $\partial_t$-invariance of $X$. Second, for a vector field $V$ we compute 

$$V(a) = \frac{1}{2a} \frac{V(g(\partial_t, \partial_t) \cdot g(X,X) - g(\partial_t, \partial_t  ) \cdot V(g(X,X))}{g^2 (X,X)},$$ 

\v this is bounded if $  V(g(\partial_t, \partial_t)) $ and $V(g(X,X))$ are bounded. In particular, for $V= \partial_t $ we have 

$$  V(g(\partial_t, \partial_t)) = 2 g (\nabla_{\partial_t} \partial_t , \partial_t) = {\rm Hess} t (\partial_t, \partial_t) , \qquad V(g(X, X)) = \dot{g}_t (X,X) < D g_t (X,X) , $$ 

\v and for $V=X $ we get

$$  X(g (\partial_t, \partial_t)) = 2 g(\nabla_X \partial_t, \partial_t) = 2 {\rm Hess} t (X, \partial_t) , \qquad X(g(X,X)) = 2 g (\nabla_X X , X) . $$

\v Taking into account the preceding equations we see that it is sufficient to show that $g(\nabla_X X, \nabla_X X) < E$ and $ g(\nabla_X X , X ) <E $ for a universal constant $E$. The way to show this will be via the metric $G:= g_+^t$ for which $ g <G $. It is thus sufficient to show $G(\nabla_X X, \nabla_X X )$ is uniformly bounded as the rest follows from the Cauchy-Schwarz inequality. On the initial hypersurface $S$, we have $\nabla_X X = 0$, in particular $G(\nabla_X X , \nabla_X X) =:u$ uniformly bounded on $S$. Let us observe the function $u$ along an integral curve of $\gr t$ and let $u'$ denote the partial derivative in the direction of $\partial_t$. We will show that $u' < F \cdot \sqrt{u} + H \cdot u + Z $ for universal constants $F,H,Z$, and then the statement follows by common ODE comparison theorems. First, by definition of $G$, we have $G(\nabla_X X , \nabla_X X) = g(\nabla_X X , \nabla_X X) + g(\nabla_X X , \partial_t)^2$ and therefore 

\bean
u' &=&  \partial_t  g ( \nabla_X X, \nabla_X X   ) +  2 g(\nabla_X X, \partial_t ) \cdot \partial_t g ( \nabla_X X , \partial_t  )\\
&=& 2 g (\nabla_t \nabla_X X, \nabla_X X) + 2  g(\nabla_X X, \partial_t ) \cdot g (\nabla_t \nabla_X X, \partial_t ) + 2 g(\nabla_X X, \partial_t ) \cdot g (\nabla_X X, \nabla_t \partial_t)
\eean

\v Let us focus on the first term. We compute

\bean
g (\nabla_t \nabla_X X, \nabla_X X) &=& g (\nabla_X \nabla_t X, \nabla_X X) + g(R(\partial_t , X) X, \nabla_X X)\\ 
&=& g (\nabla_X \nabla_X \partial_t, \nabla_X X) + g(R(\partial_t , X) X, \nabla_X X)
\eean

\v The last term of this expression, in turn, can be estimated by $I \cdot \sqrt{G(\nabla_X X, \nabla_X X)} $ for a universal constant $I$ because of the condition on the curvature. Thus it remains to be shown that the first part can be estimated in a similar manner. We write $Y:= \gr t$ and calculate

\bean
g (\nabla_X \nabla_X \partial_t, \nabla_X X)\\ 
= g(\nabla_X \nabla_X (g^{-2} ( Y , Y ) \cdot Y), \nabla_t X )\\
= g( \nabla_X ( g^{-2}(Y, Y)  \nabla_X  Y), \nabla_X X ) + g( \nabla_X (X(g^{-2} (Y , Y)  ) \cdot Y ), \nabla_X X)\\
= X(g^{-2}(Y , Y)) g( \nabla_X Y, \nabla_X X ) + g^{-2}( Y , Y) g(\nabla_X \nabla_X Y , \nabla_X X)\\
+ X(X( g^{-2}(Y, Y)    )) g(Y, \nabla_X X)  + X(g^{-2}( Y, Y )) g(\nabla_X  Y, \nabla_X X).  
\eean

\v If we denote these four additive terms as $T_1, T_2, T_3, T_4$, then we get estimates

$$\vert T_1 \vert \leq A \cdot G(\nabla_X Y, \nabla_X X)$$

\v for a universal constant $A \geq X (g^{-2} (Y,Y))$ (which exist due to the boundedness of $X$ and of the Hessian of $t$),  

$$ \vert  T_2  \vert < g^{-2} (Y,Y) \big(  \vert  G(\nabla^{(2)}_{X,X}) Y, \nabla_X X) \vert  + \vert G(\nabla_{\nabla_X X} Y, \nabla_X X) \vert \big) < B \sqrt{u} + C u $$

\v by the condition on the second derivative and the Cauchy-Schwarz inequality for the first additive term and because of the condition on the Hessian for the second one. Equally

\bean
&T_3&/g(Y, \nabla_X X)\\ 
&=& - 4 \big( g(\nabla^{(2)}_{X,X} Y, Y) g^{-3} (Y,Y)  + g(\nabla_X Y , \nabla_X Y) g^{-3} (Y,Y) - 3 g^2(\nabla_X Y, Y) g^{-4} (Y,Y) \big),
\eean

\v so that $\vert T_3 \vert$ can be estimated against $K \sqrt{u}$ for a universal constant $K$, and the bound on the last term $T_4$ follows as well from the Cauchy-Schwarz inequality and from the bound on the Hessian of $t$. Finally, bounds on the other terms in the expression for $u'$ can be obtained in complete analogy, thus the second assertion of the theorem follows. The final conclusions hold due to Theorem \ref{completelevelsets} and Theorem \ref{embed}. \hfill \qed

\section{Counterexamples}


\v The following list of simple but instructive counterexamples shows that without any further completeness assumption it might happen that no Cauchy temporal function with desired additional properties exists:

\bigskip

\v {\bf Example 1:} {\em Not every g.h. manifold $M$ admits a C time function, not even if $M$ is flat:} This elementary fact is seen by considering the causal diamond $D(p,q)$ of two causally related points $p << q $ in Minkowski space $\R^{1,d}$. The $d$-dimensional disk $S$ around the center $\frac{1}{2}(p+q)$ orthogonal to $p-q$ is a Cauchy surface for $D(p,q)$, thus any other Cauchy surface $T$ is a graph over $S$, and by the special product form of the metric there is a noncontracting diffeomorphism $\pi: T \rightarrow S$. Therefore let $\{ x_n \}_{n \in \N} $ be a nonconvergent Cauchy sequence in $S$ which exists because of noncompleteness of the disk, then the $\pi^{-1} (x_n)$ form a nonconvergent Cauchy sequence again, thus $T$ is not complete. Therefore $D(p,q)$, while globally hyperbolic, is not C. 

\bigskip

\v {\bf Example 2:} {\em Not every g.h. manifold $M$ admits an sm Cauchy temporal function, not even if $M$ is flat:} Take a Lipschitz continuous zigzag curve $Z : x \rightarrow (x, z(x)) $ in $\R^{1,1}$ as the graph of a real function $z: \R \rightarrow \R$ with $ z' \in \{ \pm 1 \}$ almost everywhere defined by $z(x) = \frac{1}{2} (x - [x])$ for $2n \leq x \leq 2n+1$, $n \in \Z$, and $z(x) = \frac{1}{2} - \frac{1}{2} (x - [x]) $ for $2n+1 \leq x \leq 2n+2$, $n \in \Z$. Then, for any smooth positive function $\phi: \R \rightarrow \R$ with $ \phi(t) = 1 $ for all $t \leq 0$, $\vert \phi' (t) \vert < \frac{1}{4}$ for all $t \in \R$ and $\lim_{t \rightarrow \infty}  \phi(t) = 0$  consider the open region $ \{ (x, y) \in \R^2 \vert  z(x) - \phi(x) < y < z(x) + \phi(x) \} $. This is a causal subset of $\R^{1,1}$ and therefore globally hyperbolic, but it is not sm: If there is a Cauchy surface given as $ C(x) := (x, c(x)) $, then it is easy to calculate that $ \vert \vert W (x) \vert \vert = c'' (x) \sqrt{\frac{1}{1 - (c' (x))^2}} = c'' (x) / \vert \vert C' (x) \vert \vert  $. Now using the intermediate value theorem, $c' $ has to change from $1/3 $ to $-1/3$ on a piece which is quenched together more and more, thus the Weingarten tensor $W$ cannot be bounded for any Cauchy surface. Now it is an easy calculation as well that $g'_t (X,Y) = 2f \cdot W(X,Y) $ and therefore $g'_t (X,Y) = 2f^2 \cdot Hess(t) (X,Y)$, $W_t (X,Y) = 2f \cdot Hess (t) (X,Y)$, for $X,Y \perp \grad (t)$, and $Hess(t) = (2f)^{-1} W_t (X,Y)$, therefore, if $f$ is bounded, the condition on the Hessian cannot be satisfied.

\bigskip

\v {\bf Example 3:} {\em Not every g.h. manifold $M$ admits an M Cauchy temporal function, not even if $M$ is flat and has a C Cauchy temporal function:} Let $\Phi: \R \rightarrow (0, \infty)$ be a smooth symmetric function (i.e., $\Phi (x) = \Phi (-x) $ for all $x \in \R$) with $\vert \Phi ' (x) \vert < \frac{1}{4}$, $\int_0^{\pm \infty} \sqrt{1 - (\Phi' (x))^2} dx = \infty$ and $\lim_{x \rightarrow \infty} \Phi (x) = 0 $. Then consider $A:= \{ (y,x) \in \R^{1,1} : - \Phi (x) < y < \Phi (x)  \}$. Clearly, $\{ 0 \} \times \R$ is a complete Cauchy surface for the flat Lorentzian manifold $A$. It can easily be seen that there are Cauchy temporal functions such that {\em all} level sets are complete: just take the Cauchy temporal function having the graphs of $t \times \Phi$, $t \in (-1,1)$, as level sets, then its length can be estimated from below by the length of the graph of $ \Phi$, which is infinite. Now, if the function $f$ in the metric decomposition was bounded globally (even if we considered a reparametrization on an interval instead of all of $\R$) there would be a universal lower bound for the length of every integral line of the gradient flow of the corresponding temporal function. But it is easily seen that for every $\e >0$ there is a point $p$ on the $x_1$ axis such that {\em every} timelike curve passing through $p$ has length smaller than $\e$.   


\bigskip

\bigskip

\v Note that in all the examples presented the incompleteness plays a fundamental role. For the last counterexample that will show that weaker assumptions than completeness do not suffice: Let a Lorentzian manifold be called {\bf maximal} if there is no nonsurjective open isometric embedding of $M$ into another Lorentzian manifold. It is an easy exercise (\cite{BEE}) that causal completeness implies maximality. The converse, however, is wrong in general.

\bigskip

\v {\bf Example 4:} {\em Not every g.h., conformally flat and maximal manifold $M$ with a complete Cauchy surface admits an M Cauchy temporal function, irrespectively of the dimension:} Let $\Phi: C^{\infty} (\R, (0, \infty))  $ with $\vert \dot{\Phi} (r) \vert < \frac{1}{2}$ for all real $r$ and with ${\rm lim}_{r \rightarrow \infty} \Phi (r) = 0 $, and let $u: (-1,1) \rightarrow [0, \infty)$ be smooth and symmetric with $u((-\frac{1}{2}, \frac{1}{2}) = \{ 0 \}$ and ${\rm lim}_{r \rightarrow \pm 1} u (r) = \infty $ but such that the integral of its square root is finite, then the conformal factor $ \psi (x_1, x_0) := 1 + u( \frac{x_0}{\Phi (x_1 )})  $ makes its finiteness region $A := \{ (x_1, x_0) \in \R^2 \vert \pm x_0 < \Phi (x_1)  \} $ (which is a causal subset of $\R^{1,1}$ and thus g.h. due to the condition $\vert \dot{\Phi} (r) \vert < \frac{1}{2}$) maximal as scalar curvature diverges along any causal geodesic. Still, the $x_1$ axis is a complete Cauchy surface as in its open neighborhood $A' := \{ (x_1, x_0) \in \R^2 \vert \pm x_0 < \frac{1}{2} \Phi (x_1)  \} $, $g$ equals the Minkowski metric. And all causal lines are of finite length due to the finite integral condition.

\bigskip

\v It is an interesting open question which geometric conditions (on completeness, curvature...) one can pose to ensure the existence of a SMC Cauchy temporal function.



\begin{thebibliography}{99}



\bibitem{BGP}
Christian B\"ar, Nicolas Ginoux, Frank Pf\"affle: {\em Wave equations on Lorentzian manifolds and quantization}. ESI Lectures in Mathematics and Physics. Z\"urich: European Mathematical Society Publishing House (2007) 

\bibitem{BEE}
John K. Beem, Paul E. Ehrlich, Kevin L. Easley: {\em Global Lorentzian Geometry}, Marcel Dekker Inc. (1996) 

\bibitem{BESSE} 
Arthur L. Besse: {\em Einstein manifolds}. Ergebnisse der Mathematik und ihrer Grenzgebiete Bd. 10. Springer-Verlag (1987). 






\bibitem{BeSa1} 
Antonio N. Bernal, Miguel S\'anchez: 
{\em On smooth {C}auchy hypersurfaces and
  {G}eroch's splitting theorem.}
\newblock { Comm. Math. Phys. 243, pp. 461-470} (2003)

\bibitem{BeSa2} 
Antonio N. Bernal, Miguel S\'anchez:
{\em Smoothness of time functions and the metric splitting of globally
hyperbolic spacetimes}, { Comm. Math. Phys. 257, pp. 43- 50} (2005)

\bibitem{BeSa} 
Antonio N. Bernal, Miguel S\'anchez: {\em Globally hyperbolic spacetimes can
be defined as ``causal" instead of ``strongly causal''}, { Class. Quant. Grav. 24, pp. 745-750} (2007)

\bibitem{BS4}
Antonio Bernal, Miguel S\'anchez: {\em Further results on the smoothability of Cauchy hypersurfaces and Cauchy time functions}, Lett.Math.Phys. 77, pp. 183-197 (2006)

\bibitem{mC}
Anna Maria Candela, Jos\'e Luis Flores, Miguel S\'anchez: {\em Global hyperbolicity and Palais-Smale condition for action functionals in stationary spacetimes}.  Adv. Math.  218 ,  no. 2, pp. 515-536  (2008). 


\bibitem{oM1} 
Olaf M\"uller: {\em The Cauchy problem of Lorentzian minimal surfaces in globally hyperbolic manifolds}. { Annals of Global Analysis and Geometry 32, no 1, pp. 67-85} (2007)


\bibitem{MS} 
Olaf M\"uller, Miguel S\'anchez: {\em Lorentz manifolds isometrically embeddable in $\mathbb{L}^n$}. Transactions of the American Mathematical Society, Vol. 363, no 10, pp. 5367-5379 (2011)


\bibitem{o5} 
Olaf M\"uller: {\em A note on closed isometric embeddings}, Journal of mathematical analysis and applications, vol. 349, no1, pp. 297-298 (2008)

\bibitem{O} 
Barrett O'Neill, {\em Semi-Riemannian Geometry with applications to Relativity}, Academic Press (1983)


\bibitem{RS} 
Alfonso Romero, Miguel S\'anchez: {\em On completeness of certain families of semi-Riemannian manifolds}, Geometriae Dedicata vol. 53, Number 1, pp. 103-117 (1994)

\bibitem{vW}
Volkmar W\"unsch: {\em Cauchy's problem and Huygens' principle for relativistic higher spin wave equations in an arbitrary curved space-time},
General Relativity and Gravitation, vol. 17, no 1, pp. 15-38 (1985)


\end{thebibliography}
\end{document}